%% file: main.tex
\pdfoutput=1

\documentclass[letterpaper, 10 pt, conference]{ieeeconf}  % Comment this line out if you need a4paper

\IEEEoverridecommandlockouts                              % This command is only needed if 
\usepackage{graphicx,algorithm,algorithmic,mathtools,enumerate,bm,hyperref,amsthm}
\usepackage{fullpage,graphicx,psfrag,amsmath,amsfonts,amssymb,pgfplots, verbatim,mathtools,dsfont,xcolor,tikz}
\usepackage[style=base]{caption}
\usetikzlibrary{calc,patterns,positioning}
\usepackage{subcaption}
\usetikzlibrary{arrows}

\usepackage{todonotes}

\allowdisplaybreaks

\input defs.tex

\title{\LARGE \bf Sensing Resource Allocation Against Data-Poisoning Attacks in Traffic Routing
}

\author{
    Yue~Yu,
    Adam~J.~Thorpe,
    Jesse~Milzman,
    David Fridovich-Keil,
    Ufuk~Topcu
\thanks{Y. Yu is with the Department of Aerospace Engineering and Mechanics at the University of Minnesota Twin Cities, Minneapolis, MN, USA. Email: {\tt yuey@umn.edu}. A. Thorpe, D. Fridovich-Keil, and U. Topcu are with the Oden Institute for Computational Engineering and Sciences, and the Department of Aerospace Engineering and Engineering Mechanics at the University of Texas at Austin, Austin, TX, USA.
Email: {\tt{adam.thorpe@austin.utexas.edu}, \tt \{dfk,utopcu\}@utexas.edu}. J. Milzman is with DEVCOM Army Research Laboratory at 2800 Powder Mill Rd, Adelphi, MD, USA. Email: {\tt{jesse.m.milzman.civ@army.mil}}. This research is jointly supported by ARO W911NF-23-1-0317, ONR N00014-22-1-2703, ONR N00014-24-S-B001, and DEVCOM ARL A2I2 CRA under grant number W911NF-23-2-0011. 
}
}

\begin{document}

\maketitle
\thispagestyle{empty}
\pagestyle{empty}

%%%%%%%%%%%%%%%%%%%%%%%%%%%%%%%%%%%%%%%%%%%%%%%%%%%%%%%%%%%%%%%%%%%%%%%%%%%%%%%%
\begin{abstract}
Data-poisoning attacks can disrupt the efficient operations of transportation systems by misdirecting traffic flows via falsified data. One challenge in countering these attacks is to reduce the uncertainties on the types of attacks, such as the distribution of their targets and intensities. We introduce a resource allocation method in transportation networks to detect and distinguish different types of attacks and facilitate efficient traffic routing. The idea is to first cluster different types of attacks based on the corresponding optimal routing strategies, then allocate sensing resources to a subset of network links to distinguish attacks from different clusters via lexicographical mixed-integer programming. We illustrate the application of the proposed method using the Anaheim network, a benchmark model in traffic routing that contains more than 400 nodes and 900 links.         
\end{abstract}

%%%%%%%%%%%%%%%%%%%%%%%%%%%%%%%%%%%%%%%%%%%%%%%%%%%%%%%%%%%%%%%%%%%%%%%%%%%%%%%%
\input{introduction/introduction}
\input{routing/routing}

\input{selection/selection}
\input{selection/optimal_routing}

\input{experiments/experiments}

\input{conclusion/conclusion}

%\input{appendix/appendix}


%%%%%%%%%%%%%%%%%%%%%%%%%%%%%%%%%%%%%%%%%%%%%%%%%%%%%%%%%%%%%%%%%%%%%%%%%%%%%%%%

%\section*{APPENDIX}
% \section*{ACKNOWLEDGMENT} Some acknowledgements.

\bibliographystyle{IEEEtran}
\bibliography{IEEEabrv,reference,data_poisoning}

\end{document}

%% file: defs.tex
\newcommand{\diag}{\mathop{\rm diag}}

\newcommand{\norm}[1]{\left\lVert#1\right\rVert}
\newcommand{\mnorm}[1]{{\left\vert\kern-0.25ex\left\vert\kern-0.25ex\left\vert #1 
    \right\vert\kern-0.25ex\right\vert\kern-0.25ex\right\vert}}

\newtheorem{assumption}{Assumption}

\newcommand{\eg}{{\it e.g.}}
\newcommand{\ie}{{\it i.e.}}

%% file: introduction/introduction.tex
\section{Introduction}

Data-poisoning attacks pose an emerging threat in intelligent transportation systems \cite{hahn2019security,mecheva2020cybersecurity}. These attacks include spoofing vehicle-level basic safety messages \cite{chen2018exposing,feng2018vulnerability}, falsifying data on users' location and points of interest used by service providers such as Google Maps \cite{eryonucu2022sybil}, or spreading misinformation via social media \cite{jamalzadeh2022protecting}. Unlike cyber-attacks that tamper with software or hardware \cite{koscher2010experimental,ghena2014green}, data-poisoning attacks manipulate driver behavior via falsifying data. They can disrupt route guidance, increase congestion, misdirect traffic towards or away from targeted areas, and lead to imbalanced use of transportation resources. 

There have been many recent studies on simulating data-poisoning attacks in transportation systems and mitigating their impacts by allocating resources for countermeasures. For example, the results in \cite{eghtesad2023hierarchical} show that multi-agent reinforcement learning can effectively find the worst-case data-injection attacks. On the other hand, Stackelberg game models can generate optimal attacks on the demand data that misdirect traffic flow towards target links in the transportation network \cite{yang2023strategic}. Finally, the results in \cite{halabi2020game} show that zero-sum games not only provide the optimal distribution of attack load with maximum impact but also resilient defense resource allocation against the worst-case attack.  

These existing results often lead to conservative routing strategies under data-poisoning attacks. They model the data-poisoning attacks as adversaries with unknown strategies. They focus on defending against the worst-case attacks subject to constraints on attack budget and intensity  \cite{eghtesad2023hierarchical,halabi2020game}. Although effective in exposing and evaluating vulnerabilities, focusing on worst-case attacks often causes conservative resource allocation \cite{halabi2020game}. A promising future direction is to reduce the unknowns in the attackers' strategies by cross-validating data using infrastructure-controlled sensors \cite{chen2018exposing}. However, to the best of our knowledge, this direction still lacks a thorough investigation in the literature.

\begin{figure}
    \centering
    \includegraphics[keepaspectratio,width=\columnwidth]{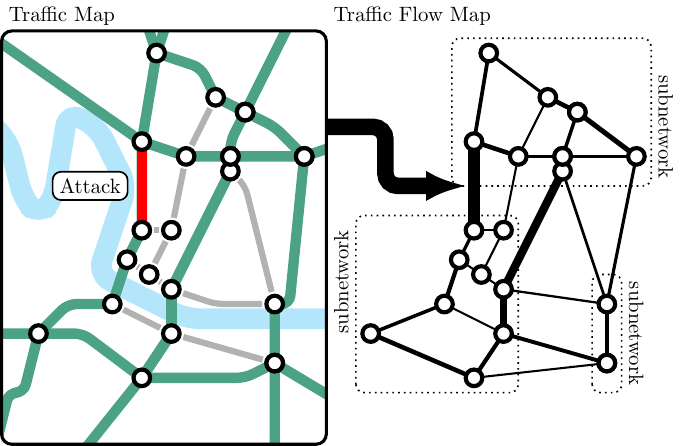}
    \caption{Data poisoning attacks can lead to incorrect flow maps, which can cause congestion or disruption in a traffic network. Our goal is to allocate sensing resources to subnetworks to detect data poisoning attacks.}
    \label{fig: head}
\end{figure}

We formulate a novel sensing resource allocation problem to facilitate traffic routing under data-poisoning attacks. The idea is to allocate limited sensing resources, \eg, the limited flight time of surveillance drones, to validate the traffic levels on a subset of links in a transportation network when the reported traffic levels are under data-poisoning attacks. Instead of completely unknown attacks, we consider attacks sampled from one of finitely many hypothetical distributions. Different distributions affect different network links with different intensities. The goal is to first allocate limited sensors to a subset of links, then use sensor observations to reduce the uncertainties on the type of attacks, and finally route traffic flows between given origin-destination pairs subject to the reduced uncertainties on the type of attacks.   

We propose a novel method to solve the resource allocation problem based on \emph{attack clustering} and mixed-integer programming. First, we compute the optimal routing under the attacks sampled from each hypothetical distribution. Second, we cluster these distributions based on the corresponding routing, such that attack distributions leading to similar routing belong to the same cluster. Finally, we propose a mixed-integer programming approach to allocate sensing resources to subnetworks that optimally distinguish attacks sampled from distinct clusters. We demonstrate the application of the proposed algorithm using the Anaheim network, a benchmark model in traffic routing that contains more than 400 nodes and 900 links.       

%% file: routing/routing.tex
\section{Traffic routing and data-poisoning attacks}
The system-optimal network routing problem is an optimization problem where a traffic planner minimizes the total travel cost of the traffic flow created by network users. This optimization problem is composed of three components: the network structure, the user demands between different origin-destination pairs, and the link cost function that captures how traffic volume affects travel time.

\subsection{Transportation network model}

We model the transportation network using a directed transportation graph \(\mathcal{G}\) composed of \(n_n\in\mathbb{N}\) nodes, each of which corresponds to the intersection of roads, and \(n_l\in\mathbb{N}\) links, each of which corresponds to a road segment. Each link is an ordered pair of distinct nodes, where the first and second nodes are the ``tail" and ``head" of the link, respectively.

An origin-destination (OD) pair is a pair of distinct nodes \((i, j)\). We let \(n_d\) denote the total number of OD pairs. Furthermore, we let vector \(d\in\mathbb{R}^{n_d}_{\geq 0}\) denote the demand of these pairs such that its \(i\)-th entry \(d_i\) is the amount of demand for the \(i\)-th OD pair. 

A \emph{route} is a sequence of links where each link's head is connected to the next link's tail. We let \(n_r\) denote the total number of routes in the network. Throughout, we only consider routes that connect the given OD pairs with positive demand, such that each route connect one unique OD pair.   

\subsection{Incidence matrices}

To formulate the network routing problem, we model the network structure using the link-route incidence matrix and the route-OD incidence matrix.

\subsubsection{Link-route incidence matrix} We describe the link-route incidence relation in graph \(\mathcal{G}\) using the notion of \emph{link-route incidence matrix}. In particular, we let matrix \(F\in\mathbb{R}^{n_l\times n_r}\) be such that its entry \(F_{ij}\) associates link \(i\) with route \(j\) as follows:
\begin{equation}\label{eqn: lr incidence}
     F_{ij}=\begin{cases}
    1, & \text{if link \(i\) is on route \(j\),}\\
    0, & \text{otherwise.}
    \end{cases}
\end{equation}

\subsubsection{Route-OD incidence matrix} We describe the route-OD incidence relation in graph \(\mathcal{G}\) using the notion of \emph{route-OD incidence matrix}. In particular, we define a matrix \(H\in\mathbb{R}^{n_d\times n_r}\) be such that the entry \(F_{ij}\) associates the \(i\)-th OD pair with route \(j\) as follows:
\begin{equation}\label{eqn: rod incidence}
     H_{ij}=\begin{cases}
    1, & \text{if route \(j\) connects the \(i\)-th OD pair,}\\
    0, & \text{otherwise.}
    \end{cases}
\end{equation}

\subsection{Link cost function}
The link cost function is a continuously differentiable function that describes how the cost of using a link increases with the traffic volume on the link. One of the most commonly used link cost functions is the Bureau of Public Roads (BPR) model \cite{patriksson2015traffic}. In this model, we suppose that the traffic volume on link \(j\) is given by \(f_j+y_j\), where \(f_j\in\mathbb{R}_{\geq 0}\) is the ambient traffic volume, independent of the traffic planner's decision, and \(y_j\in\mathbb{R}_{\geq 0}\) is the traffic volume assigned to link \(j\) by the traffic planner. Then the cost of using link \(j\), denoted by \(\phi_j(f_j+y_j)\), is defined as:
\begin{equation}
\label{eq: cost kernel}
\phi_j(f_j+y_j)\coloneqq b_j+w_j\left(\frac{f_j+y_j}{c_j}\right)^4,
\end{equation}
where \(b_j\in\mathbb{R}_{\geq}\) denotes the \emph{nominal travel cost} when the traffic volume on link \(j\) is zero, \(c_j\in\mathbb{R}_{\geq 0}\) denotes the \emph{nominal capacity} of link \(j\), \(w_j\in\mathbb{R}_{>0}\) is the weighting parameter that accounts for link congestion.

\subsection{System-optimal traffic routing}

System-optimal routing aims to find a link flow assignment that serves the user demand for each OD pair, satisfies flow conservation constraints at each node, and minimizes the total link cost. To formulate this routing problem as an optimization problem, we introduce the following variables:
\begin{equation}
    \label{eqn: flow variables}
     y\in\mathbb{R}^{n_l}, \enskip z\in\mathbb{R}^{n_r}.
\end{equation}
In particular, vector \(y\) denotes the \emph{link flow vector} whose \(i\)-th entry \(y_i\) denotes the total amount of flow on link \(i\). Vector \(z\) denotes the \emph{route flow vector} whose \(j\)-th entry \(z_j\) denotes the amount of flow on route \(j\).

Using \eqref{eqn: flow variables}, we formulate the system-optimal routing problem as the following optimization 
\begin{equation}
\label{opt: routing}
\begin{array}{ll}
\underset{y, z}{\mbox{minimize}} &  \sum_{j=1}^{n_l} y_j\phi_j(f_j+y_j)\\
\mbox{subject to} & Hz=d, \, Fz=y,\, z\geq 0_{n_r}.
\end{array}
\end{equation}
The objective function in \eqref{opt: routing} measures the total cost of link flows. The constraints in \eqref{opt: routing} ensure flow conservation at each node.

\subsection{Data-poisoning attacks against routing}

We consider the scenario where the traffic planner cannot observe the accurate ambient traffic flow \(f\in\mathbb{R}^{n_l}_{\geq 0}\) before assigning traffic flow via solving optimization~\eqref{opt: routing}. Instead, it only observes an \emph{attacked ambient flow}, denoted by 
\begin{equation}
\hat{f}\coloneqq f+a,
\end{equation}
where \(a\in\mathbb{R}^{n_l}\) denotes an additive attack vector. Such an attack is due to false reports on traffic conditions, which typically originate from Sybil-based attackers that aim to divert traffic away or towards the targeted area \cite{eryonucu2022sybil}. 

If the traffic planner simply computes an optimal flow assignment by solving optimization~\eqref{opt: routing} with \(f\) replaced by \(\hat{f}\), the attack vector can dramatically affect the resulting flow assignment.

%% file: selection/selection.tex
\section{Sensing resource allocation against data-poisoning attacks}

We now formulate a sensor allocation problem where the traffic planner allocates sensors to observe the effects of the data-poisoning attacks on a subset of links. We assume that the attack comes from one of many hypothetical distributions. The goal of sensor allocation is to pinpoint which hypothesis is most likely. To this end, we start with the following assumption on the distributions of the attacks. For simplicity, we assume the attacks on different links are independent. 
\begin{assumption}\label{asp: Gaussian attack}
   The traffic planner knows that there exists \(\mu^i\in\mathbb{R}^{n_l}\) and \(\sigma^i\in\mathbb{R}_{>0}^{n_l}\) for \(i=1, 2, \ldots, n_a\) such that, for any \(1\leq j\leq n_l\), \(\mu_j^i\leq \hat{f}_j\) for any \(1\leq i\leq n_a\). Furthermore,  
    \(a_j\sim\mathcal{N}(\mu^i_j, (\sigma^i_j)^2)\) for some \(1\leq i\leq n_a\).
\end{assumption}
We also make the following assumptions on the subnetworks where the traffic planner can deploy sensors and the corresponding cost. 
\begin{assumption}\label{asp: subnetworks}
    The traffic planner can allocate sensors to \(n_g\in\mathbb{N}\)  distinct subnetworks to detect the attacks on these subnetworks. Each subnetwork contains a subset of links, denoted by \(\mathcal{L}_1, \mathcal{L}_2, \ldots, \mathcal{L}_{n_g}\subset \{1, 2, \ldots, n_l\}\), where \(|\mathcal{L}_i|=n_s^i\) for all \(i=1, 2, \ldots, n_g\).
\end{assumption}

\begin{assumption}\label{asp: resources}
   Allocating sensors to the \(i\)-th subnetwork costs \(c_i\in\mathbb{R}_{>0}\) amount of resources. The traffic planner has a budget of \(\gamma\in\mathbb{R}_{>0}\) amount of resources.
\end{assumption}

Under Assumption~\ref{asp: subnetworks} and Assumption~\ref{asp: resources}, we aim to answer the following question: how should the traffic planner allocate its resource budget by deploying sensors to the optimal set of subnetworks? 

To answer this question, we introduce the concept of \emph{attack clustering}, which groups the types of attacks according to the traffic planner's optimal responses to attacks. Second, we define a \emph{difference function} that evaluates each subnetwork in terms of distinguishing different attacks in different clusters. Finally, we formulate a mixed-integer program to optimally allocate sensing resources to subnetworks to maximize the differences within the allowed budget.

\subsection{Best response routing and the clustering of attacks}
We start with the basic case where the traffic planner knows the type of attack. Suppose that Assumption~\ref{asp: Gaussian attack} holds and the planner knows the attack type \(i\). Then the planner's best response to the attack is to optimize the network flow by minimizing the expected value of the objective function in \eqref{opt: routing} where \(f=\hat{f}-a\) with \(a\sim\mathcal{N}(\mu^i, \diag(\sigma^i\odot \sigma^i))\), where \(\odot\) denotes the Hadamard product. In this case, we say that \(\hat{z}^i\in\mathbb{R}^{n_r}_{\geq 0}\) is the \emph{best response flow for type-\(i\) attacks} if there exists \(\hat{y}^i\in\mathbb{R}^{n_l}_{\geq 0}\) such that \((\hat{y}^i, \hat{z}^i)\) is an optimal solution of the following optimization problem
\begin{equation}
\label{opt: best response}
\begin{array}{ll}
\underset{y, z}{\mbox{minimize}} &  \sum_{j=1}^{n_l} y_j\mathds{E}[\phi_j(y_j+\hat{f}_j-a_j)]\\
\mbox{subject to} &  Hz=d, \, Fz=y,\, z\geq 0_{n_r},\\
& a\sim\mathcal{N}(\mu^i, \diag(\sigma^i)).
\end{array}
\end{equation}

Since \(a_j\sim\mathcal{N}(\mu^i_j, (\sigma^i_j)^2)\) for any \(1\leq j\leq n_l\), we know that \(y_j+\hat{f}_j-a_j^i\sim\mathcal{N}(y_j+\hat{f}_j-\mu_j^i, (\sigma_j^i)^2)\). Therefore we can show that
\begin{equation}\label{eqn: exp link cost}
\begin{aligned}
   \psi_j^i(y_j)\coloneqq &\mathds{E}[\phi(y_j+\hat{f}_j-a_j)] \\
   =&b_j+\frac{w_j}{(c_j)^4}(3(\sigma_j^i)^4+(y_j+\hat{f}_j-\mu_j^i)^4)\\
    &+\frac{6w_j}{(c_j)^4}(y_j+\hat{f}_j-\mu_j^i)^2(\sigma_j^i)^2.
\end{aligned} 
\end{equation}

Hence optimizing~\eqref{opt: best response} is equivalent to the following
\begin{equation}
\label{opt: best response 2}
\begin{array}{ll}
\underset{y, z}{\mbox{minimize}} &  \sum_{j=1}^{n_l} y_j\psi_j^i(y_j)\\
\mbox{subject to} &  Hz=d, \, Fz=y,\, z\geq 0_{n_r}.
\end{array}
\end{equation}

\input{figs/BPmap}

\begin{comment}

    {
\color{blue}
Rather than take the estimate above, we may instead compute the best optimal response directly.

Taking $f = \hat{f} - a$, and assuming attack type $t=i$, we have the following:
\begin{align}
    \nonumber & \mathbb{E} \left[ \sum_{j=1}^{n_l} \phi_j(y_j ; f_j) \; \bigg| \; t=i \right] \\
    % \nonumber =& \sum_{j=1}^{n_l}\mathbb{E} \left[ \phi_j (y_j ; f_j) \; \bigg| \; t=i \right] \\
    % \nonumber =& \sum_{j=1}^{n_l} b_j y_j + \frac{y_j}{c_j^4} \; \mathbb{E} \left[ \left( a^i_j - (\hat{f}_j + y_j) \right)^4 \right] \\
    % \nonumber =& \sum_{j=1}^{n_l} b_j y_j + \frac{y_j}{c_j^4} \left[  \sum_{k=0}^2 {4 \choose 2k } (y_j - \tilde{\mu}_j^i)^{2k} (\sigma_j^i)^{4-2k} \right] \\
    % =& \sum_{j=1}^{n_l} b_j y_j + \frac{y_j}{c_j^4} \left[  (y_j - \tilde{\mu}_j^i)^4 + 6 (\sigma_j^i)^2 (y_j - \tilde{\mu}_j^i)^2 + (\sigma_j^i)^4 \right] \\
    \label{eq:best_response_flow_polynomial}
    =& \sum_{j=1}^{n_l} \sum_{k=1}^5 \zeta^i_{j,k} \,  y_j^k
    % \intertext{where}
    % & \alpha_{j,k} = \delta_{1,k} b + c_j^{-4}\left( \right)
\end{align}
where the coefficients $\zeta^i_{j,k} = \zeta^i_{j,k}(\mu_j^i, \sigma_j^i, \hat{f}_j, c_j, b_j)$ are functions of given in Appendix~\ref{app:coefficients}.

}
%%%%%%%%
\end{comment}
%%%%%%%% Jesse's insert

Notice that different types of attacks could lead to similar best response flows. For example, if all the OD pairs concentrate in the southern part of a transportation network, then different types of attacks localized at a link in the northern part of the network will likely have similar best response flows. In other words, depending on the demand, different attack types may form \emph{clusters}. Within each cluster, different attacks may lead to similar best response flows. We illustrate this phenomenon in Fig.~\ref{fig: flow clustering}. 

To identify attack types that lead to similar best response flows, we propose a \emph{attack clustering problem}. Let \(\hat{z}^i_\star\) denote the best response flow for type-\(i\) attacks. We compute \(n_c\in\mathbb{N}\) clusters of best responses, denoted by \(v^1, v^2, \ldots, v^{n_c}\) by solving the following optimization problem:
\begin{equation}\label{opt: cluster}
    \begin{array}{ll} 
    \underset{v^1, v^2, \ldots, v^{n_c}}{\mbox{minimize}}  &  \sum_{j=1}^{n_c}\sum_{i=1}^{n_a}  \norm{v^i-\hat{z}^j_\star}_1
    \end{array}
\end{equation}
where \(n_c\in\mathbb{N}\) is a parameter denoting the total number of clusters. Here we choose the \(\ell_1\) norm because it measures the sum of the total variation of the route flow. By choosing an appropriate value for \(n_c\), the solution of the optimization problem in~\eqref{opt: cluster} ensures that, for any given \(\epsilon\in\mathbb{R}_{>0}\) and \(1\leq i\leq n_a\), we have \(\norm{\hat{z}^i-v^j}_1\) is below a certain threshold for some \(1\leq j\leq n_c\) and all \(1\leq i\leq n_a\). A trivial choice of \(n_c\) is to let \(n_c=n_a\). In practice, we search for the smallest number of clusters \(n_c\) such that each \(\hat{z}^i\) is sufficiently close to the center of at least one cluster. In the following, we let 
\begin{equation}\label{eqn: cluster i}
    \mathcal{C}_j\coloneqq \{i| \norm{\hat{z}^i-v^j}_1\leq \epsilon, 1\leq i\leq n_a \} 
\end{equation}
for all \(1\leq j\leq n_c\). Throughout, we assume that \(\epsilon\) in \eqref{eqn: cluster i} is sufficiently small such that \(\mathcal{C}_i\cap \mathcal{C}_j=\emptyset\) for all \(1\leq i, j\leq n_c\) with \(i\neq j\).

\subsection{Difference function}

Equipped with the clustering of attacks, we are ready to define a \emph{difference function} on the subnetworks. This function evaluates each subnetwork based on how well it differentiates the attacks in different clusters. 

To this end, we let \(\mathcal{P}\) denote the set of attack type pairs, where the two types in each pair belong to different clusters. We define \(\mathcal{P}\) as follows
\begin{equation}
    \mathcal{P}\coloneqq  \{(i, j)| 1\leq i<j\leq n_a, i\neq j\} \setminus \mathcal{Q}
\end{equation}
where the set \(\mathcal{Q}\) contains all of the pair of attack types from the same cluster, \ie, 
\begin{equation}
    \mathcal{Q}\coloneqq \{(i, j)| \exists 1\leq k\leq n_c, \text{ s.t. } i\in\mathcal{C}_k, j\in\mathcal{C}_k \}.
\end{equation}

In addition, we introduce a \emph{selection matrix} \(S^i\in\{0, 1\}^{n_s^i\times n_l}\) for each \(i=1, 2, \ldots, n_g\). Entry \(S^i_{kj}\) associates set \(\mathcal{L}_i\) and link \(j\) as follows:
\begin{equation}\label{eqn: selection matrix}
     S_{kj}^i=\begin{cases}
    1, & \text{if link \(j\) is the \(k\)-th link in \(\mathcal{L}^i\),}\\
    0, & \text{otherwise.}
    \end{cases}
\end{equation} 
Under Assumption~\ref{asp: Gaussian attack}, we know that 
\begin{equation}
    S^j a\sim\mathcal{N}(\xi^{ij}, \Lambda^{ij} ) 
\end{equation}
for all \(i=1, 2, \ldots, n_a\) and \(j=1, 2, \ldots, n_g\), where
\begin{equation}
    \xi^{ij}\coloneqq S^j\mu^i, \enskip \Lambda^{ij}\coloneqq S^j \diag(\sigma^i) (S^j)^\top. 
\end{equation}

We now introduce a \emph{difference matrix}, denoted by \(M\in\mathbb{R}^{n_p\times n_g}\) where \(n_p\coloneqq |\mathcal{P}|\). The entry \(M_{ij}\) associates the \(i\)-th pair in the set \(\mathcal{P}\) and subnetwork \(j\) as follows
\begin{equation}\label{eqn: metric matrix}
     M_{ij}= (\xi^{pj}-\xi^{qj})^\top(\Lambda^{pj}+\Lambda^{qj})^{\dagger}(\xi^{pj}-\xi^{qj})
\end{equation}
where \((p, q)\) is the \(i\)-th pair in the set of \(\mathcal{P}\). The value \(M_{ij}\) gives the Kullback–Leibler divergence between distribution \(\mathcal{N}(\xi^{pj}, \Lambda^{pj} )\) and distribution \(\mathcal{N}(\xi^{qj}, \Lambda^{qj} )\). Intuitively, it measures the difference of attack type \(p\) and \(q\) when only the links in subnetwork \(\mathcal{L}_j\) are observable.

\subsection{Optimal sensor allocation problem}

We propose a resource allocation problem to select the optimal subnetworks that help distinguish among different attack clusters. To this end, we introduce a binary vector of variables, denoted by
\begin{equation}
    x\in\{0, 1\}^{n_g}
\end{equation}
such that \(x_i=1\) if and only if sensors are deployed to the \(i\)-th subnetwork. We propose to compute the value of variable \(x\) using the following mixed-integer program
\begin{equation}\label{opt: sensor selection}
    \begin{array}{ll}
    \underset{x, u}{\mbox{maximize}} &  f(u) \\
    \mbox{subject to} & Mx =u,\enskip q^\top x\leq \gamma, \enskip x\in\{0, 1\}^{n_g}. 
    \end{array}
\end{equation}
Here \(f:\mathbb{R}^{n_p}\to\mathbb{R}\) is a function that evaluates the distinguishability of all the pairs of attacks in set \(\mathcal{P}\), and \(q\in\mathbb{R}^{n_p}\) is a vector whose \(i\)-th entry denotes the cost of allocating to the \(i\)-th subnetwork. Two intuitive choices for the function \(f\) are the average and the elementwise minimum function, given as follows
\begin{subequations}
    \begin{align}
        f_{\text{avg}}(u) & =\frac{1}{n_p} \mathbf{1}_{n_p}^\top u, \label{eqn: avg}\\
        f_{\text{min}}(u) & =\min_{1\leq i\leq n_p} u_i.\label{eqn: max-min}
    \end{align}
\end{subequations}
Here \( f_{\text{avg}}\) evaluates on average how distinguishable are the attack pairs in \(\mathcal{P}\), whereas \( f_{\text{min}}\) evaluates the worst-case distinguishability among all attack pairs in \(\mathcal{P}\). 

In practice, there often exist multiple allocations that maximize either \(f_{\text{avg}}\) or \(f_{\text{min}}\). In addition, to ensure resilience against the worst-case scenario, we prioritize \(f_{\text{min}}\) over \(f_{\text{avg}}\). Therefore, we propose a \emph{lexicographic approach}. In particular, we compute the optimal allocation \(x\) by solving the following mixed-integer programming:
\begin{equation}\label{opt: allocation upper}
    \begin{array}{ll}
    \underset{x, u}{\mbox{maximize}} &  f_{\text{avg}}(u) \\
    \mbox{subject to} & Mx =u,\enskip q^\top x\leq \gamma, \,u\geq \alpha \mathbf{1}_{n_p},\\
    & x\in\{0, 1\}^{n_g}. 
    \end{array}
\end{equation}
where \(\alpha\) is the optimal value of the following optimization:
\begin{equation}\label{opt: allocation lower}
    \begin{array}{ll}
    \underset{\tilde{x}, \tilde{u}}{\mbox{maximize}} &  f_{\text{min}}(\tilde{u}) \\
    \mbox{subject to} & M\tilde{x} =\tilde{u},\enskip q^\top \tilde{x}\leq \gamma, \enskip \tilde{x}\in\{0, 1\}^{n_g}. 
    \end{array}
\end{equation}
The idea is to first narrow down the search space to the optimizers of the problem in \eqref{opt: allocation lower}, then search for the allocation that further solves optimization~\eqref{opt: allocation upper}.  

%% file: figs/BPmap.tex
   \tikzset{
    node/.style={rectangle,draw,minimum size=1.5em},link/.style={->}
}
%,
 
\begin{figure}
    \centering
    \begin{tikzpicture}
    \coordinate (f1) at (4.5, 2.45);
    \coordinate (f3) at (4.5, 1.95);
    \coordinate (f2) at (4.5, 0.25);
    \coordinate (f4) at (4.5, -0.25);

   \node[label={\footnotesize{Attack Types}}] at (0, 3.4) (L1) {}; 
   \node[label={[align=center]\footnotesize{Best Response Flow}}] at (4.7, 3.4) (L1) {}; 
   
   \node[node] at (0, 3) (A1) {\footnotesize Type 1}; 
   \node[node] at (0, 2) (A2) {\footnotesize Type 2}; 
   \node[node] at (0, 1) (A3) {\footnotesize Type 3}; 
   \node[] at (0, 0.2) (D) {\(\vdots\)}; 
   \node[node] at (0, -0.8) (A4) {\footnotesize Type \(n_a\)}; 

   %\draw [black,fill=gray!5] (4.5,1) ellipse (2cm and 2.6cm);
   \draw [blue,dashed] (4.7,2.2) ellipse (0.8cm and 0.8cm);
   \draw [red,dashed] (4.7,0) ellipse (0.8cm and 0.8cm);

   \node[] at (4.7, 1.2) (C1) {\footnotesize \color{blue} Cluster 1}; 
    \node[] at (4.7, -1) (C1) {\footnotesize \color{red} Cluster 2}; 

   \fill[blue] (f1) circle (1.5pt);
   \node[label=right:{$z^1$}] at (f1) {};
   \fill[blue] (f3) circle (1.5pt);
   \node[label=right:{$z^2$}] at (f2) {};
   \fill[red] (f2) circle (1.5pt);
   \node[label=right:{$z^3$}] at (f3) {};
   \fill[red] (f4) circle (1.5pt);
   \node[label=right:{$z^{n_a}$}] at (f4) {};
   
   %\draw [black,dashed] (5,1) ellipse (2cm and 2.5cm);
   
   \draw[gray, dashed, ->] (A1) -- (f1);
   \draw[gray, dashed, ->] (A2) -- (f2);
   \draw[gray, dashed, ->] (A3) -- (f3);
   \draw[gray, dashed, ->] (A4) -- (f4);
   % \draw[link] (2) -- (3);

    \end{tikzpicture} 
    \caption{Different types of attacks can lead to clusters of best response flows}
    \label{fig: flow clustering}
\end{figure}
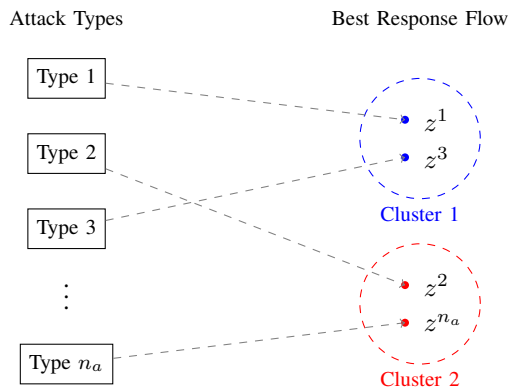

%% file: selection/optimal_routing.tex
\subsection{Optimal routing with sensing observations}

After allocating the sensing resources to subnetworks, we can collect observations of the ground-truth attacks on these subnetworks and use them to obtain an approximation of the ground-truth routing objective function in \eqref{opt: routing}. In particular, we let \(x^\star\) denote an optimal solution of the optimization in \eqref{opt: allocation upper} and \(\mathcal{M}\subset \mathcal{L}\) be the corresponding set of links that got an allocation of sensing resources, \ie,  
\begin{equation}
    \mathcal{M}\coloneqq \{k|\exists \mathcal{L}_j \text{ such that } k\in\mathcal{L}_j \text{and } x_j^\star>0 \}. 
\end{equation}
In addition, we let \(n_m\coloneqq |\mathcal{M}|\) and introduce the \emph{sensing selection matrix} \(E\in\mathbb{R}^{n_m\times n_l}\) whose \(ij\)-th entry is as follows:
\begin{equation}
    E_{ij} = \begin{cases}
    1, & \text{if \(j\in\mathcal{M}\),}\\
    0, & \text{otherwise.}
    \end{cases}
\end{equation}
Let \(a\in\mathbb{R}^{n_l}\) denote the ground truth attack on all links and 
\begin{equation}
    o\coloneqq Ea
\end{equation}
denote the vector of attacks observed in \(\mathcal{M}\).
If \(a\sim\mathcal{N}(\mu^i, \diag(\sigma^i\odot \sigma^i)\), then 
\begin{equation}
    o\sim\mathcal{N}(E\mu^i, E\diag(\sigma^i\odot \sigma^i)E^\top).
\end{equation}

Given the above observations, we approximate the optimization in \eqref{opt: routing} as follows:
\begin{equation}
\label{opt: routing after sensing}
\begin{array}{ll}
\underset{y, z}{\mbox{minimize}} &  \sum_{j\in\mathcal{M}} y_j\phi_j(y_j+\hat{f}_j-a_j)\\
&+\sum_{j\in\mathcal{L}\setminus\mathcal{M}} \sum_{i=1}^{n_a} \omega^iy_j\psi_j^i(y_j)\\
\mbox{subject to} &  Hz=d, \, Fz=y,\, z\geq 0_{n_r},
\end{array}
\end{equation}
where \(\sum_{i=1}\omega^i =1\) and
\begin{equation}\label{eqn: likelihood}
\omega^i>0,\,\omega^i\propto \frac{\exp\left(-\frac{1}{2}\norm{E\mu^i-o}_{E\diag(\sigma^i\odot \sigma^i)E^\top}^2\right)}{\sqrt{(2\pi)^{n_m}\det(E\diag(\sigma^i\odot \sigma^i)E^\top)}},
\end{equation}
for all \(1\leq i\leq n_a\), where \(\odot\) denotes the Hadamard product. The central idea of this approximation is as follows. For each link \(j\), if \(j\in\mathcal{M}\), we can observe the corresponding attack \(a_j\) and use the same link cost function as the one in \eqref{opt: routing}. If \(j\notin\mathcal{M}\), then we approximate the link cost function using a weighted sum of the expected link cost used in \eqref{opt: best response}. The \(i\)-th summand is the expected link cost under type-\(i\) attack, first introduced in \eqref{eqn: exp link cost}. We weight the  \(i\)-th summand with \(\omega^i\), which, according to \eqref{eqn: likelihood}, is proportional to the likelihood of observing attack \(o\) in link set \(\mathcal{M}\) if \(a\sim\mathcal{N}(\mu^i, \diag(\sigma^i\odot \sigma^i))\).

%% file: experiments/experiments.tex
\section{Numerical experiments}
We demonstrate the application of the proposed sensor allocation algorithm using the \emph{Anaheim network}, a benchmark network model in traffic routing \cite{transpnet} based on the city of Anaheim, California. In the following, we discuss the details of the network model, attack type clustering, and the corresponding sensing allocation.

\begin{figure}[!htb]
    \centering
    \includegraphics[trim={5cm 5cm 5cm 4cm},clip,width=2.1in]{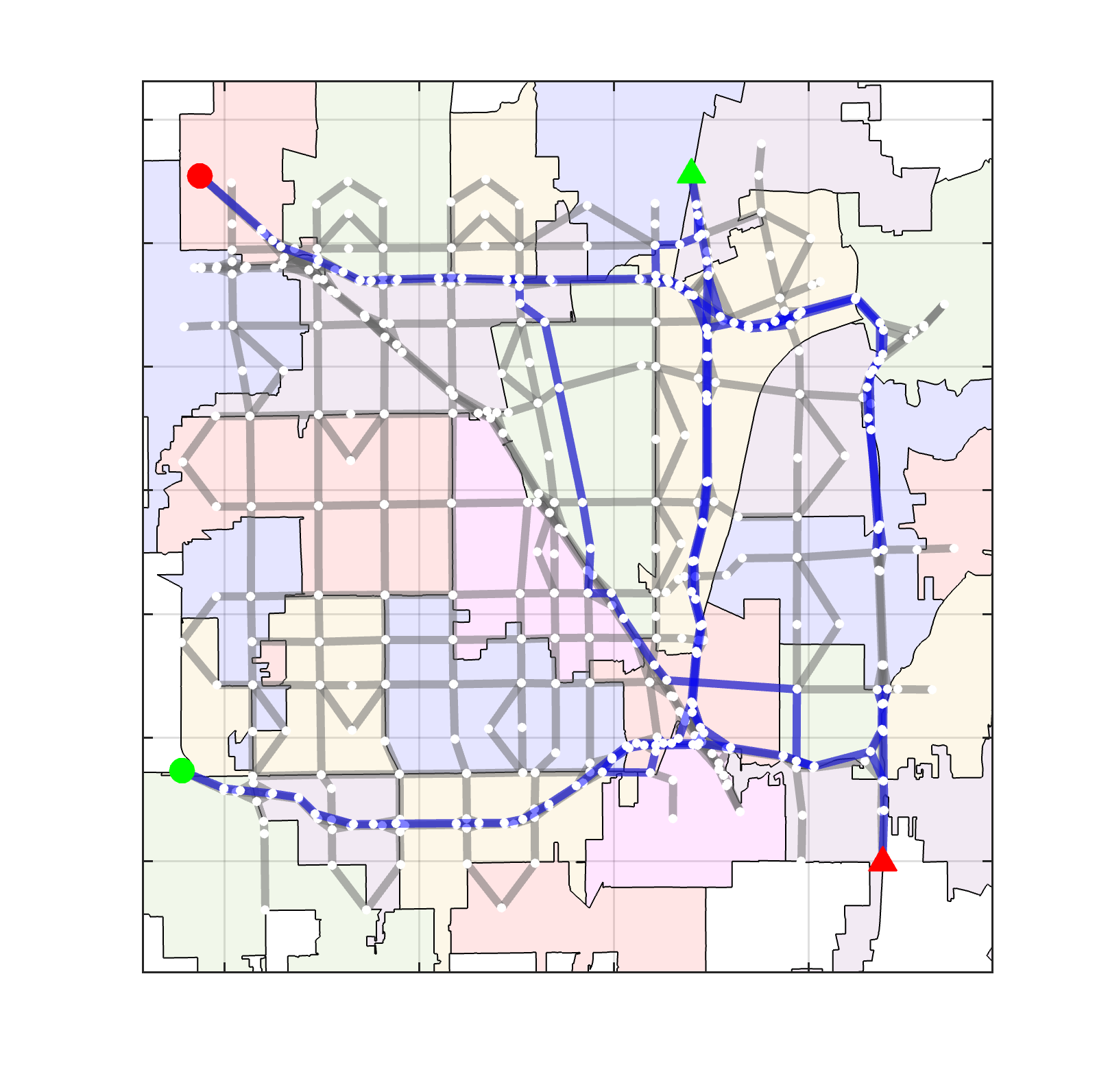}
  \caption{An illustration of the Anaheim network model. The blue links are on the routes between two origins (red and green circles) and two destinations (red and green triangles) are highlighted. }\label{fig: map}
\end{figure}

\begin{figure*}[!htb]
    \centering
   \begin{subfigure}{2in}
  \includegraphics[trim={5cm 3cm 5cm 3cm},clip,width=\linewidth]{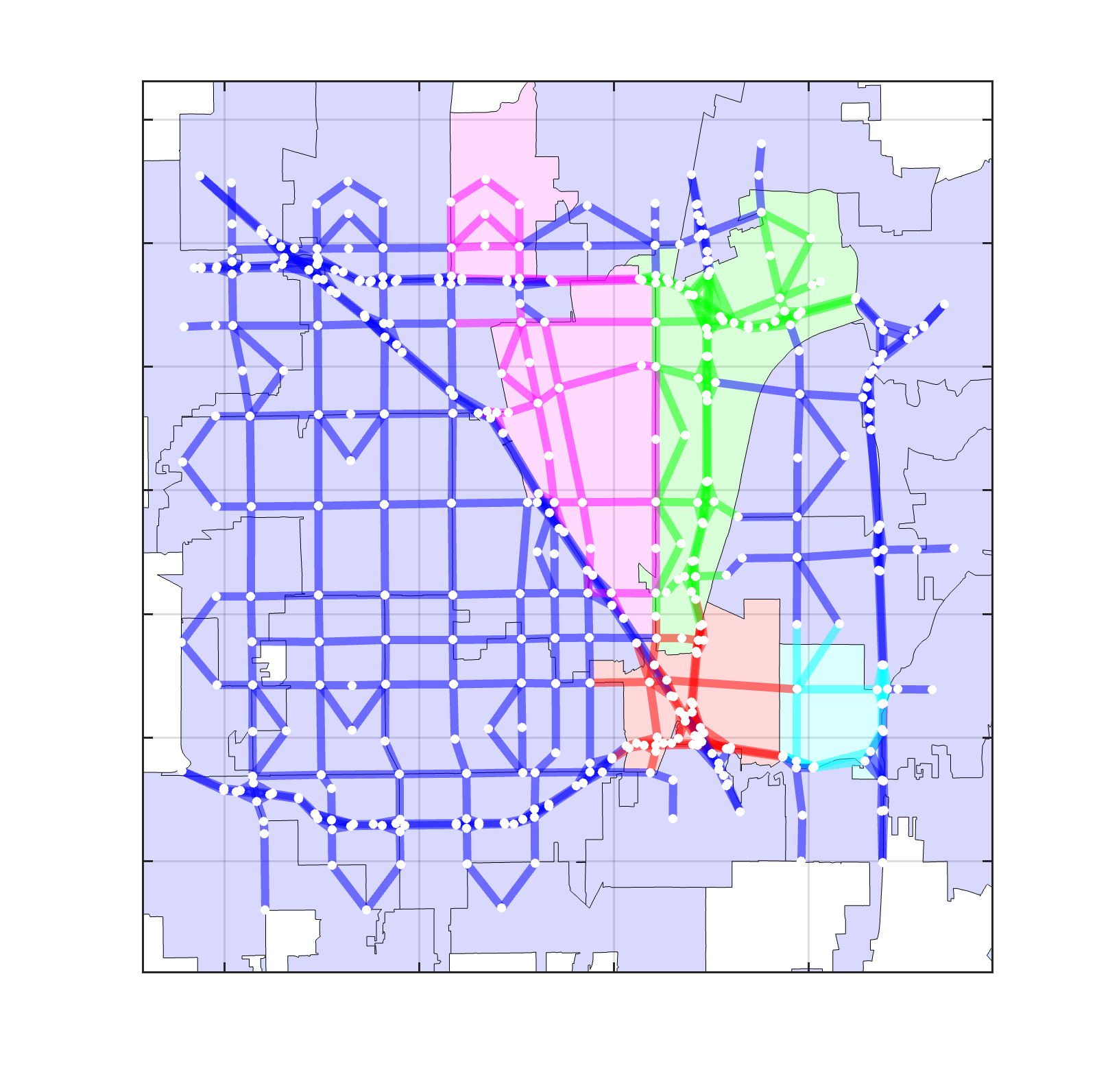}
  \caption{High budget (\(\gamma=27\)).}\label{fig: max budget}
  \end{subfigure}
  \hfill
  \begin{subfigure}{2in}
  \includegraphics[trim={5cm 3cm 5cm 3cm},clip,width=\textwidth]{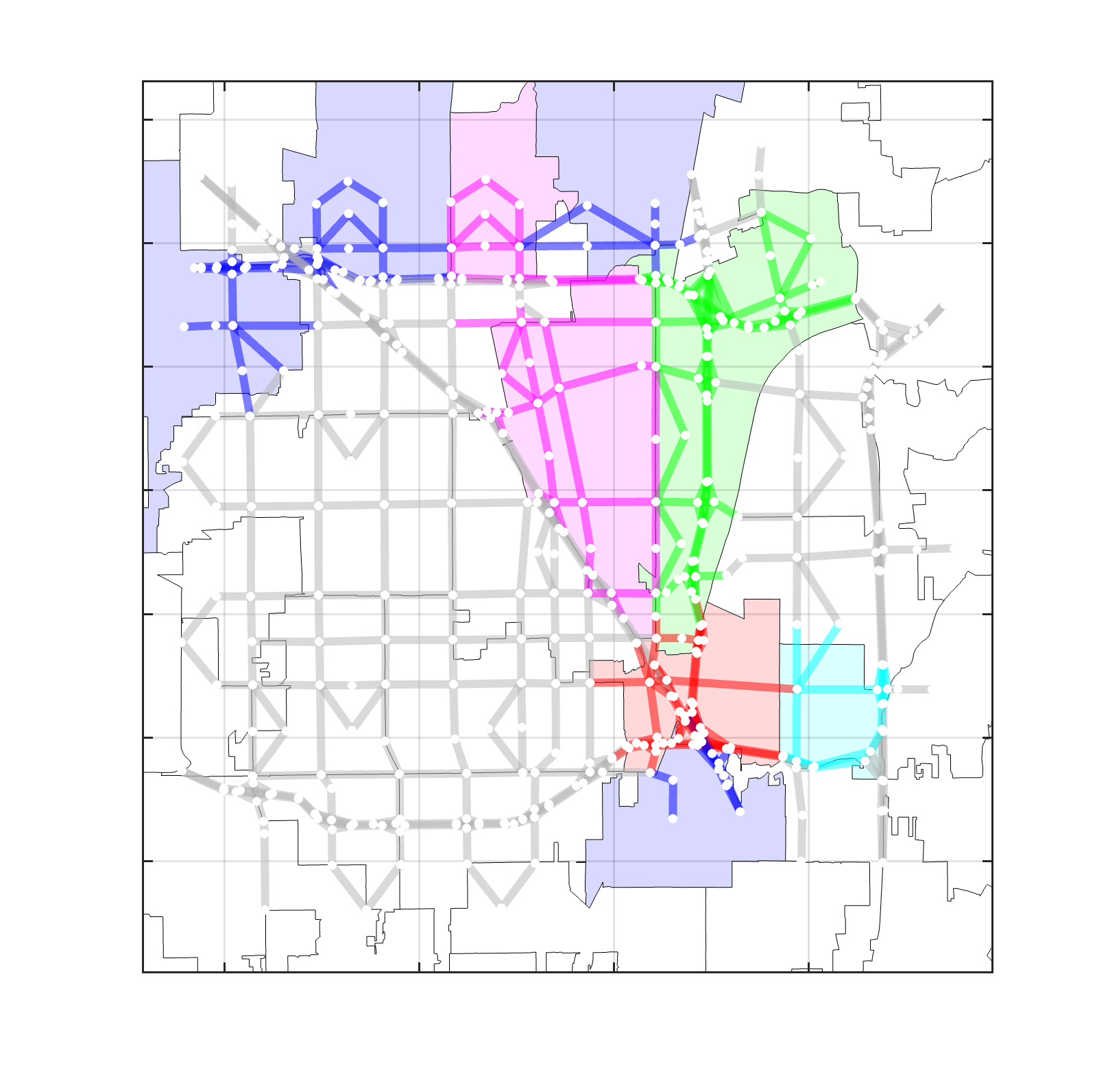}
  \caption{Medium budget (\(\gamma=9\)).} \label{fig: mid budget}
  \end{subfigure} 
 \hfill
\begin{subfigure}{2in}
  \includegraphics[trim={5cm 3cm 5cm 3cm},clip,width=\linewidth]{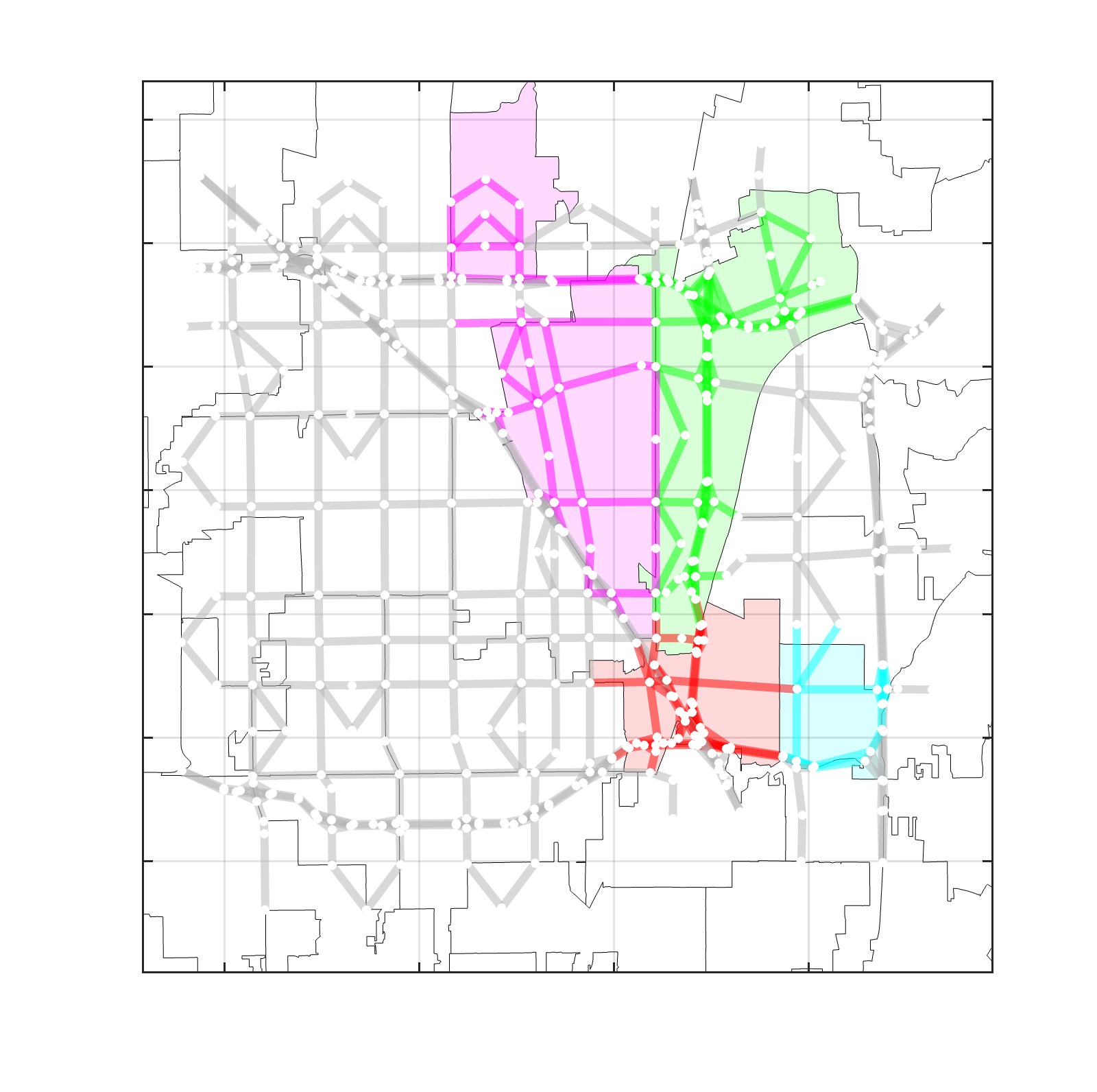}
  \caption{Low budget (\(\gamma=5\)).}\label{fig: min budget}
  \end{subfigure}
  
   \caption{Sensor allocation to detect data poisoning attacks in unknown zip code areas under different budgets \(\gamma\), which appears in \eqref{opt: allocation upper} and \eqref{opt: allocation lower}. We highlight the subnetworks selected by~\eqref{opt: allocation upper}, where different colors indicate attacks in different clusters.}
   \label{fig: alloaction}
\end{figure*}

\subsection{Anaheim network and the attack model}

The Anaheim network is a benchmark network model in the traffic routing literature \cite{transpnet}. The network contains \(n_n=416\) nodes and \(n_l=914\) directed links. Within this network, we consider a total of \(n_r=8\) routes between two distinct origin-destination pairs. See Fig.~\ref{fig: map} for an illustration. The model also includes the nominal travel cost \(b\in\mathbb{R}_{\geq 0}^{n_l}\) and nominal capacity \(c\in\mathbb{R}_{\geq 0}^{n_l}\). See \cite{transpnet} for the details of these parameters. Throughout we assume that the ground-truth value of the ambient traffic volume is half of the nominal capacity, \ie, \(f=\frac{1}{2} c\).

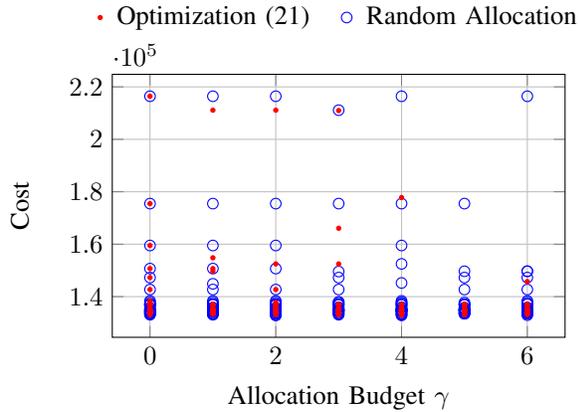
\begin{figure}[!htb]
    \centering
    \begin{tikzpicture}
        \begin{scope}
    
            \begin{axis}[
                grid=major,
                ylabel=Cost,
                xlabel=Allocation Budget \(\gamma\),
                height=2in,
                width=0.95\columnwidth,
                scatter/classes={
                    opt={mark=*,red,mark size=0.8pt},
                    rnd={mark=o,blue}
                },
                legend style={draw=none,at={(0,0)},anchor=north west},
                legend columns=2,
                /tikz/column 2/.style={
                    column sep=10pt,
                },
            ]
            \addplot[
                scatter,
                only marks,
                scatter src=explicit symbolic,
            ]
            table[x=y,y=x,meta=label] {./experiments/data/new_data.dat};
            \legend{~Optimization \eqref{opt: allocation upper},~Random Allocation}
            \end{axis}
    
        \end{scope}
    \end{tikzpicture}
    \caption{The distribution of the optimal value of optimization~\eqref{opt: routing after sensing} under different types of attacks. For each budget, we plot the cost values corresponding to \(n_a=27\) different types of attacks, one for each zip code area in Fig.~\ref{fig: map}.}
    \label{fig: cost}
\end{figure}

Based on this network, we construct both the types of attack as well as the subnetwork partition for sensor allocation based on the zip code areas. The Anaheim network covers 27 different zip codes (see Fig.~\ref{fig: map} for an illustration). We partition the network into \(n_g=27\) subnetworks according to the zip code of each node and link. Furthermore, we consider \(n_a=27\) different types of attacks. For each type \(i\) where \(1\leq i\leq n_a\),  we let \(N^i\in\mathbb{R}^{n_l\times n_l}\) such that
\begin{equation}
     N_{ij}^i=\begin{cases}
    1, & \text{if \(i=j\) and \(i\in\mathcal{L}^i\),}\\
    0, & \text{otherwise.}
    \end{cases}
\end{equation} 
and Assumption~\ref{asp: Gaussian attack} holds with
\begin{equation}\label{eqn: attack gaussian}
    \mu^i=30 N^ic, \enskip \sigma^i=(1/10) c.
\end{equation}
for all \(1\leq i\leq n_a\). In other words, the type \(i\)-attack will only affect the links in the \(i\)-th subnetwork---which are the links in the \(i\)-th zip code--- by falsifying a false increase in traffic volume sampled from the Gaussian distribution who means and variances are given by \eqref{eqn: attack gaussian}.  

\subsection{Attack clustering and sensor allocation}
We cluster the \(n_a=27\) different types of attacks by solving optimization~\eqref{opt: cluster} and determine which zip code areas to allocate sensors by solving optimization~\eqref{opt: allocation upper}. We use the commercial solver Mosek \cite{aps2019mosek} with default tolerance setting to solve optimization~\eqref{opt: cluster}.  We illustrate these results using Fig.~\ref{fig: alloaction}. In particular, Fig~\ref{fig: max budget} shows that, when given enough budget, the solution of optimization~\eqref{opt: allocation upper} suggests allocating sensing resources to all subnetworks. Fig~\ref{fig: max budget} also shows that the \(n_a=27\) different type of attacks form \(n_c=5\) clusters. Because the candidate routes only use a small subset of links (see Fig.~\ref{fig: map} for an illustration), attacks in most areas do not lead to different best response flows, and hence belong to one large cluster (painted blue). On the other hand, links in seven zip code areas (painted with non-blue colors) cause changes in the best response flows. A possible explanation is that links in these areas appear in different candidate routes.

Fig.~\ref{fig: alloaction} also shows the priorities of the proposed allocation. In particular, Fig.~\ref{fig: mid budget} and Fig.~\ref{fig: min budget} show that as we decrease the budget value \(\gamma\), the allocation suggests areas that correspond to different clusters rather than those from the same clusters. This agrees with the idea of optimization~\eqref{opt: allocation upper}, which is to distinguish attacks in different clusters.  

Finally, we illustrate the optimal values of optimization~\eqref{opt: routing after sensing} under different allocation budgets and different types of attacks in Fig.~\ref{fig: cost}. Furthermore, we compare these values against the case where the areas suggested by optimization~\eqref{opt: allocation upper} are replaced by random ones. Fig.~\ref{fig: cost} shows that even though there are \(n_g=27\) different subnetworks, sensing five of them is sufficient to reduce the routing cost to a close-to-optimal value. On the other hand, the random allocation can lead to a much higher cost, depending on which type of attack is in effect.

%% file: conclusion/conclusion.tex
\section{Conclusion}
We introduce a resource allocation method for traffic routing under data-poisoning attacks with unknown types. The idea is to first cluster different types of attacks according to the corresponding best-response flows, then allocate sensing resources to distinguish attacks in different clusters via lexicographical mixed integer programming.

However, the current work still has limitations. For example, we only consider Gaussian-type attack distributions that are independent across different links in Assumption~\ref{asp: Gaussian attack}. This assumption does not account for multimodal attack distributions or correlated attacks on multiple links. In addition, it remains unclear how the computation of the proposed method scales with network size. Our future work aims to address this limitation, as well as extend the current results to resource allocation for dynamic routing in stochastic networks.